\documentclass[a4paper,10pt,openright]{article}
\usepackage[T1]{fontenc}
\usepackage[applemac]{inputenc}
\usepackage{mathrsfs}
\usepackage{indentfirst}
\usepackage{amssymb,amsmath}
\usepackage{amsthm}
\usepackage{makeidx}
\usepackage{array}
\usepackage{graphicx}
\usepackage{latexsym}
\usepackage{multirow}
\usepackage{caption}
\usepackage{fancyhdr}
\usepackage{bm}
\usepackage{mathptmx}      
\usepackage[T1]{fontenc}
\usepackage[applemac]{inputenc}
\usepackage{mathrsfs}
\usepackage{indentfirst}
\usepackage{amssymb,amsmath}
\usepackage{amssymb}
\usepackage{amsfonts}
\usepackage{amscd}
\usepackage{enumerate}
\usepackage{amsthm}
\usepackage{makeidx}
\usepackage{array}
\usepackage{graphicx}
\usepackage{latexsym}
\usepackage{caption}
\usepackage{fancyhdr}
\usepackage{bm}

\newcommand{\beq}{\begin{equation}}
\newcommand{\eeq}{\end{equation}}
\newcommand{\bea}{\begin{eqnarray}}
\newcommand{\eea}{\end{eqnarray}}
\newcommand{\Prob}{\mathbb{P}}

\newcommand{\Natural}{\mathbb{N}}

\def\<{\langle}
\def\>{\rangle}
\newtheorem{proposition}{Proposition}
\newtheorem{rem}{Remark}

\begin{document}

\title{Bivariate semi-Markov Process for Counterparty Credit Risk}

\author{
Guglielmo D'Amico\\  
        Universit\'a di Chieti-Pescara G. D'Annunzio,\\
        Dipartimento di Farmacia\\
        g.damico@unich.it\\
\and
Raimondo Manca\\  
        Sapienza Universit\'a di Roma\\
        Dipartimento di Metodi e Modelli per l'Economia, il Territorio e la Finanza\\
         raimondo.manca@uniroma1.it\\
\and
Giovanni Salvi\\  
        Sapienza Universit\'a di Roma\\
        Dipartimento di Metodi e Modelli per l'Economia, il Territorio e la Finanza\\
         giovanni.salvi@uniroma1.it\\}

\maketitle

\begin{abstract}
\indent We consider the problem of constructing an appropriate multivariate model for the study of the counterparty credit risk in credit rating migration problem. For this financial problem different multivariate Markov chain models were proposed. However the markovian assumption may be inappropriate for the study of the dynamic of credit ratings which typically show non markovian like behaviour.  
In this paper we develop a semi-Markov approach to the study of the counterparty credit risk by defining a new multivariate semi-Markov chain model. Methods are given for computing the transition probabilities, reliability functions and the price of a risky Credit Default Swap.
\end{abstract}

{\textbf{keywords:}} counterparty credit risk; multivariate semi-Markov chains; algorithm.

\section{Introduction}
\label{intro}
The current financial crisis has stressed the importance of the study of the correlations in the financial market. In this regard, the study of the risk of default of the counterparty, in any financial contract, has become crucial in the credit risk. For a complete treatment about credit risk we refer to the classical book of Bielecki and Rutkowski \cite{Bielecki}. Many works have been done trying to describe the counterparty risk in a Credit Default Swap (CDS) contract, but all these works are based on the Markovian approach to the credit risk, see for example Crepey et al.[6]  or Ching and Ng \cite{ChingBook}.\\
\indent It has been shown that, the Markov chain based models, are too restrictive for the description of accurate rating dynamics, see for example Carty and Fons \cite{Cart94}. Indeed, they require that the distribution functions of the sojourn times in a rating class before to have a transition should be exponentially or geometrically distributed for continuous and discrete time models, respectively.\\
\indent  In an attempt to produce more efficient credit rating models semi-Markov processes were proposed for the first time as applied to credit ratings in the paper by \cite{DamicoManca0}; more recent results were given in \cite{DamicoManca4,DamicoManca2}.\\
\indent It is important to dispose of efficient migration models because reliable rating prediction is of interest for pricing rating sensitive derivatives \cite{Vasi06}, \cite{DamicoManca3}, for the valuation of portfolio of defaulting bonds, for credit risk management and capital allocation.\\
\indent No results are available for counterparty credit risk for semi-Markov credit rating models. Such an extension is not straightforward as randomness in the sojourn times and memory effects are to be appropriately managed.\\
\indent In this paper, therefore, we present a novel multivariate semi-Markov model to approach the counterparty risk in a CDS contract. In Section 2 we investigate the behaviour of the multivariate process in the transient case and we derive equations for the transition probabilities with backward recurrence time dependence and reliability functions. An algorithm useful to perform the computations of the transition probabilities is provided. In Section 4 we analyze the counterparty credit risk in a CDS contract. Section 5 is dedicated to a summary of this study and suggestions for future outlooks.

\section{Multivariate Semi-Markov Chains} \label{Sec:SMC}

\begin{figure}
\begin{center}
\includegraphics[scale=0.65]{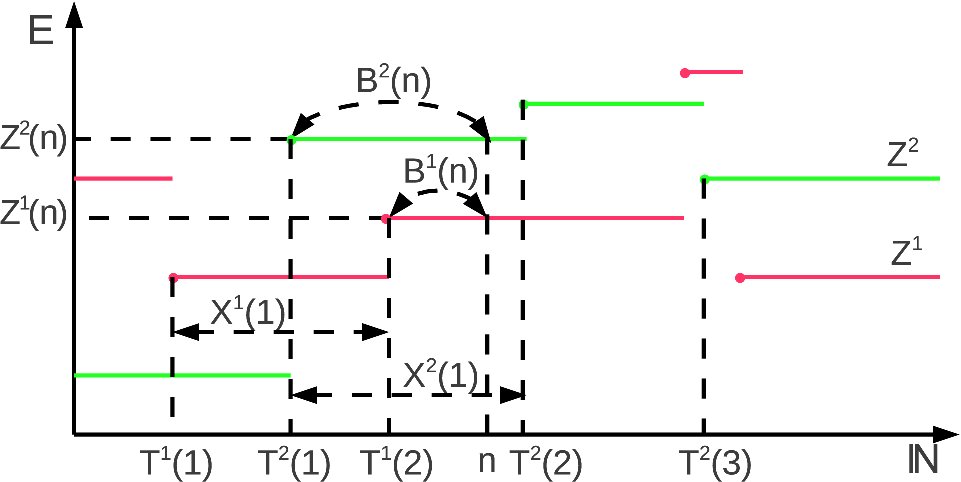}
\caption{The trajectory of the double-component system is shown as a function of time. In the picture sojourn times, transition times and backward recurrence times are shown.}
\label{Fig:trj2SMP}
\end{center}
\end{figure}

The main issue of this session is to define Multivariate Semi-Markov Chains (MVSMC), for this purpose we first introduce Markov Renewal chains, for more details see for example Barbu and Limnios \cite{Barbu}, Janssen and Manca \cite{Manca} and Limnios and Opri\c an \cite{Limnios}.\\
\indent Let $J = (J_{n})_{n \in \Natural}$ be a Markov Chain with values in $E = \{ 1, \dots ,d \}$, that is
$$\Prob (J_{n+1}=j \ | \ J_{0}=i_{0},\dots,J_{n}=i) = \Prob (J_{n+1}=j \ | \ J_{n}=i) = p_{ij} \qquad \forall \  i,j \in E.$$
The matrix $P = (p_{ij})_{i,j \in E}$ is the one step transition probability matrix. Suppose that the permanence in the state is triggered by a renewal moment process $T = (T_{n})_{n \in \Natural}$, with values in $\Natural$, defined by $T_{n+1} = T_{n} + X_{n} \qquad \forall \  n \in \Natural.$\\
\indent $X = (X_{n})_{n \in \Natural}$ is a sequence of random variables representing the sojourn times in the $n^{th}$-state. They have a conditional cumulative distribution function given by
$$\Prob (X_{n+1} \leq t \ | \ J_{n}=i, J_{n+1}=j) = F_{ij}(t) \qquad \forall \  i \in E.$$
The couple $(J,T) = (J_{n},T_{n})_{n \in \Natural}$ is said to be a Markov Renewal Chain (MRC) if
$$\Prob (J_{n+1}=j,T_{n+1}-T_{n} \leq k \ | \ (J_{0},T_{0}),\dots,(J_{n},T_{n})) = \Prob (J_{n+1}=j,X_{n+1} \leq k \ | \ J_{n}) .$$
\indent In this case, $J = (J_{n})_{n \in \Natural}$ is said to be the associated embedded Markov Chain. \\ 
\indent Further let's define the counting process $N(t) = \max \{ n \in \Natural \ | \ T_{n} \leq t \},$ which gives the number of transitions of the MRC up to time t. \\ 
\indent 
Now we are ready to introduce the Semi-Markov Chain. Let $Z = (Z_{n})_{n \in \Natural}$ be the process, with value in $E = \{ 1,...,d \}$, defined by
\begin{eqnarray}
Z_{k} := J_{N(k)}, \qquad k \in \Natural .
\label{Def:SMC}
\end{eqnarray}
\indent The process (\ref{Def:SMC}) is called Semi-Markov Chain (SMC) associated with the MRC $(J,T)$, whose cumulated kernel is denoted by $\mathbf{Q} = (Q(k) ; k \in \Natural)$ and defined for all $i,j \in E$ and $k \in \Natural$, by
$$ Q_{ij}(k) := \Prob (J_{n+1}=j,X_{n+1} \leq k \ | \ J_{n}=i). $$
\indent The element $Q_{ij}(k)$ is the probability that the system makes next transition in state $j$ with sojourn time less or equal to k given that the present state is $i$.\\
\indent Now, we are ready to introduce the multivariate semi-Markov chain. Let us consider a system consisting of $\gamma$ parts, each part has values in $E = \{ 1, \dots, d \}$. 
Let's call 
$J^{\alpha} = (J^{\alpha}_{n})_{n \in \Natural}$, for $\alpha = 1, \dots ,\gamma$, the sequence of states visited by the $\alpha$ part with values in $E$. 
Denote with $(T^{(\alpha)}_{n})_{n \in \Natural}$ the sequence of transition times of the $\alpha$-th component with state space $\Natural$. Let also introduce the sequence of random variables $X^{\alpha}_{n} = T^{\alpha}_{n+1} - T^{\alpha}_{n}$, for every $n \in \Natural$. $X^{\alpha}_{n}$ is the sojourn time in state $J^{\alpha}_{n}$.\\
\indent
Let us define also the counting processes
$$N^{\alpha}(k) = \max \{ n \in \Natural \ | \ T^{\alpha}_{n} \leq k \} \qquad \forall \alpha = 1, \dots , \gamma \ \textrm{and} \ k \in \Natural ,$$
which give us the number of transitions of the part $\alpha$ up to time $k$; we denote by $\mathbf{N}(k) = (N^{1}(k), \dots , N^{\gamma}(k))$, the vector of all such numbers. \\
\indent
To define the MVSMC we make two assumptions named in the following A1 and A2.
\begin{itemize}
\item[A1] The process $\mathbf{J} = (J^{\alpha})_{\alpha = 1, \cdots ,\gamma}$ is a multivariate Markov chain in the sense described below. \\ 
\indent Given the component $\alpha$, we denote the vector of all components except of $\alpha$ by the symbol $-\mathbf{\alpha}=(1,2,...,\alpha -1, \alpha +1,...,\gamma)$. We introduce the marginal one step transition probability for the multivariate Markov chain $\mathbf{J}$ for all $s \in \Natural$ as
\begin{eqnarray}
\label{Eq:TransMatrixMSMC}
& \Prob ( J_{N^{\alpha}(s)+1}^{\alpha} = j \ | \ \sigma(J_{h}^{\alpha} ,h\leq N^{\alpha}(s)), \sigma(J_{h}^{-\mathbf{\alpha}} ,h\leq \mathbf{N}^{-\alpha}(s)+1)) \nonumber \\
& =\Prob ( J_{N^{\alpha}(s)+1}^{\alpha} = j \ | \ J_{N^{\alpha}(s)}^{\alpha} , J_{\mathbf{N}^{-\mathbf{\alpha}}(s)}^{-\mathbf{\alpha}})=:p^{\alpha}_{\mathbf{J}_{\mathbf{N}(s)},j}(s)
\end{eqnarray}
where $\sigma(J_{h}^{\alpha} ,h\leq N^{\alpha}(s))$ denotes the natural filtration of $J^{\alpha}$ and $\sigma(\mathbf{J}^{-\mathbf{\alpha}} ,h\leq \mathbf{N}^{-\alpha}(s)+1))$ is the natural filtration of the multivariate $\mathbf{J}^{-\mathbf{\alpha}}$ process. We assume also that the process is time homogeneous and then, the transition probabilities do not depend on time $s$, so we have
\begin{eqnarray}
p^{\alpha}_{\mathbf{i},j}(s) = p^{\alpha}_{\mathbf{i},j} \qquad \forall \ s \in \Natural .
\end{eqnarray}
\item[A2] For every $\alpha , \beta = 1, \cdots ,\gamma$, the sequences of sojourn times $(X^{\alpha}_{n})_{n \in \Natural}$ and $(X^{\beta}_{n})_{n \in \Natural}$ are independent in the sense of formula (\ref{Eq:IndepX}) here below: \\
\end{itemize}
\begin{eqnarray}
F^{\alpha}_{i_{\alpha}}(k;s) & := & \Prob (X_{N^{\alpha}(s)+1}^{\alpha} \leq k \ | \ \sigma(J_{N^{\beta}(h)}^{\beta},X_{N^{\beta}(h)}^{\beta}),h \leq s ,\beta = 1,\dots ,\gamma) \nonumber \\
& = & \Prob (X_{N^{\alpha}(s)+1}^{\alpha} \leq k \ | \ J_{N^{\alpha}(s)}^{\alpha} = i_{\alpha})   \nonumber \\ 
& = & \Prob (X_{n+1}^{\alpha} \leq k \ | \ J_{n}^{\alpha} = i_{\alpha}) = F^{\alpha}_{i_{\alpha}}(k) \qquad \forall \ s \in \Natural,
\label{Eq:IndepX}
\end{eqnarray}
\noindent
where $\sigma(J_{N^{\beta}(h)}^{\beta},X_{N^{\beta}(h)}^{\beta}, h \leq s, \beta = 1,\dots ,\gamma)$ is the natural filtration of the multidimensional $(\mathbf{J},\mathbf{X})$ process. \\
\indent
Then, in our model we suppose that the sequence of sojourn times $(X^{(\alpha)}_{n})_{n \in \Natural}$ depends only on the visited state of the $\alpha$-component $J_{n}^{\alpha}$.
In the second step we use the independence of the $\alpha$ part sojourn time from the trajectory of the other component of the system.\\
\indent The discrete time semi-Markov kernel for each part $\alpha$ of the system can be formally defined as follow
\begin{eqnarray}
Q_{\mathbf{i},j}^{\alpha} (k;s) := \Prob ( J_{N^{\alpha}(s)+1}^{\alpha} = j, X_{N^{\alpha}(s)+1}^{\alpha} \leq k \ | \ J_{N^{1}(s)}^{1} = i_{1}, \dots , J_{N^{\gamma}(s)}^{\gamma} = i_{\gamma}) ,
\label{Eq:DefKernel}
\end{eqnarray}
for all $s \in \Natural$, and where $\mathbf{i} = (i_{1}, \dots , i_{\gamma})$  is a vector in $E^{\gamma}$ and $j$ is an element of $E$. By hypothesis we have that the semi-Markov kernel does not depend on time $s$, that is
\begin{eqnarray}
Q_{\mathbf{i},j}^{\alpha} (k;s) = Q_{\mathbf{i},j}^{\alpha} (k) \qquad \forall \ s \in \Natural .
\end{eqnarray}
\indent Notice that the semi-Markov kernel introduced in formula (\ref{Eq:DefKernel}) can be written as
\begin{eqnarray}
Q_{\mathbf{i},j}^{\alpha} (k) = p^{\alpha}_{\mathbf{i},j} F^{\alpha}_{i_{\alpha}}(k) ,
\end{eqnarray}
where it is implied that $\mathbf{i} = \mathbf{J}_{\mathbf{N}(s)+1}$ and $j = J_{N^{\alpha}(s)+1}$, for all $s \in \Natural$. \\
\indent
In what follow, it can be useful to speak in terms of probability to have a transition exactly at a certain time, hence we define
\begin{eqnarray}
q_{\mathbf{i},j}^{\alpha} (k;s) &  =  & \Prob ( J_{N^{\alpha}(s)+1}^{\alpha} = j, X_{N^{\alpha}(s)+1}^{\alpha} = k \ | \ J_{N^{1}(s)}^{1} = i_{1}, \dots , J_{N^{\gamma}(s)}^{\gamma} = i_{\gamma})  \nonumber \\
& = & Q_{\mathbf{i},j}^{\alpha} (k + 1) - Q_{\mathbf{i},j}^{\alpha} (k) =: q_{\mathbf{i},j}^{\alpha} (k) ,
\label{Eq:qMSMC}
\end{eqnarray}
for all $s \in \Natural$. Notice that the transition probabilities (\ref{Eq:TransMatrixMSMC}) can be expressed in terms of $(q^{\alpha}_{\mathbf{i},j})$ as
\begin{eqnarray}
p^{\alpha}_{\mathbf{i},j} = \sum_{k = 0}^{\infty} q_{\mathbf{i},j}^{\alpha} (k) .
\end{eqnarray}
\indent We define, for each component, for every state $i \in E$ and time $k \in \Natural$, the unconditional sojourn time cdf:
\begin{eqnarray}
H_{i}^{\alpha}(k) & := & \Prob (X^{\alpha}_{n+1} \leq k \ | \ J_{n}^{\alpha} = i) = \Prob (X_{N^{\alpha}(s)+1}^{\alpha} \leq k \ | \ J_{N^{\alpha}(s)}^{\alpha} = i_{\alpha}) \nonumber \\
& = & \sum_{j \in E} Q_{ij}^{\alpha}(k) \qquad \forall \ s \in \Natural ,
\end{eqnarray}
\indent We notice that in this case we have $H^{\alpha} = F^{\alpha}$, that is
\begin{eqnarray}
H_{i}^{\alpha}(k) = F_{i}^{\alpha}(k) = \sum_{j \in E} Q_{ij}^{\alpha}(k) .
\label{Eq:HequalF}
\end{eqnarray}
\indent The evolution of the multivariate $(\mathbf{J}, \mathbf{X})$ process can be described as follow, given an initial state the next state occupied by the system is determined according to the evolution of the multivariate Markov chain while the sojourn time in the present state,  of every part of the system, is determined according to the joint distribution of $\mathbf{X}$.
In the bivariate case an example of trajectory is shown in Figure (\ref{Fig:trj2SMP}).\\
\indent Now, we can introduce the MVSMC associated to the kernel $(Q^{\alpha})_{\alpha = 1, \cdots , \gamma}$. We denote the MVSMC by $\mathbf{Z} = (Z^{\alpha})_{\alpha = 1, \dots , \gamma}$ where  $\forall \alpha =1,\ldots,\gamma$
\begin{eqnarray}
Z^{\alpha}(k) := J^{\alpha}_{N^{\alpha}(k)} \qquad \forall k \in \Natural .
\label{Eq:DefMSMC}
\end{eqnarray}
\indent
The transition function of the MSMC for $\alpha$-th component is defined by
\begin{eqnarray}
\Phi_{\mathbf{i},j}^{\alpha}(k) = \Prob ( Z^{\alpha}(k) = j \ | \ \mathbf{Z}(0) = \mathbf{i}) ,
\label{Eq:DefPhiMSMC}
\end{eqnarray}
for all $\mathbf{i} \in E^{\gamma}$, $j \in E$ and $k \in \Natural$. \\
\indent
In the semi-Markov environment, the Markovian property is preserved only at transition times, in general we have to take into account for the history of the process. Then transition probabilities depend on the elapsed time in a given state, see for example D'Amico et al. \cite{DamicoManca}. \\
\indent
In general when we look for transition functions of semi-Markov model we need to introduce the concept of backward recurrence time. \\
\indent
The backward recurrence time for any component $\alpha$, denoted by $B^{\alpha}_{t}$, specifies at any time the age of the state of the $\alpha$-th part, that is
$$B^{\alpha}_{t} := t - T^{\alpha}_{N^{\alpha}(t)} \qquad \forall \alpha = 1, \dots ,\gamma .$$
\indent In other words, $B^{\alpha}_{t}$ gives the time from the last transition of $\alpha$-th component. \\
\indent
The one step transition probabilities for component $\alpha$, by tacking into account for bacward recurrence time, is defined by
\begin{eqnarray}
\Prob (J^{\alpha}_{N^{\alpha}(s)+1} = j, X^{\alpha}_{N^{\alpha}(s)+1} = k + v_{\alpha} \ | \ \mathbf{J}_{\mathbf{N}(s)} = \mathbf{i} \ , \ \mathbf{T}_{\mathbf{N}(s)} = s - \mathbf{v} \ , \ \mathbf{T}_{\mathbf{N}(s) + \mathbf{1}} > s ) =: q^{\alpha}_{\mathbf{i},j} (v_{\alpha},k) \nonumber
\end{eqnarray}
where, $v_{\alpha}$ stand for the backward of component $\alpha$ at time s, $\mathbf{i} = (i_{1},\dots,i_{\gamma})$ and $\mathbf{J}_{\mathbf{N}(\cdot)} = (J^{1}_{N_{1}(\cdot)},\dots,J^{\gamma}_{N_{\gamma}(\cdot)})$ are vectors in $E^{\gamma}$, $\mathbf{v} = (v_{1},\dots,v_{\gamma})$ and $\mathbf{T}_{\mathbf{N}(\cdot)} = (T^{1}_{N_{1}(\cdot)},\dots,T^{\gamma}_{N_{\gamma}(\cdot)})$ are vectors in $\Natural^{\gamma}$. \\
\indent
\begin{proposition}
The one step transition probabilities for component $\alpha$ are given by
\begin{eqnarray}
q^{\alpha}_{\mathbf{i},j} (v_{\alpha},k) = \frac{F^{\alpha}_{i_{\alpha}}(k+v_{\alpha}) - F^{\alpha}_{i_{\alpha}}(k+v_{\alpha} - 1)}{1 - H^{\alpha}_{i_{\alpha}}(v_{\alpha})} \cdot p_{\mathbf{i},j}^{\alpha} \ .
\label{Def:qBack}
\end{eqnarray}
\end{proposition}
\begin{proof}
By applying Bayes rules we get
\begin{eqnarray}
& & \Prob (J^{\alpha}_{N^{\alpha}(s)+1} = j, X^{\alpha}_{N^{\alpha}(s)+1} = k + v_{\alpha} \ | \ \mathbf{J}_{\mathbf{N}(s)} = \mathbf{i} \ , \ \mathbf{T}_{\mathbf{N}(s)} = s - \mathbf{v} \ , \ \mathbf{T}_{\mathbf{N}(s) + \mathbf{1}} > s )  \nonumber \\
& = & \Prob (X^{\alpha}_{N^{\alpha}(s)+1} = k + v_{\alpha} \ | \ J^{\alpha}_{N^{\alpha}(s)+1} = j \ , \ \mathbf{J}_{\mathbf{N}(s)} = \mathbf{i} \ , \ \mathbf{T}_{\mathbf{N}(s)} = s - \mathbf{v} \ , \ \mathbf{T}_{\mathbf{N}(s) + \mathbf{1}} > s )  \nonumber \\
& \times & \Prob (J^{\alpha}_{N^{\alpha}(s)+1} = j \ | \ \mathbf{J}_{\mathbf{N}(s)} = \mathbf{i} \ , \ \mathbf{T}_{\mathbf{N}(s)} = s - \mathbf{v} \ , \ \mathbf{T}_{\mathbf{N}(s) + \mathbf{1}} > s ) .
\end{eqnarray}
Then, by using assumptions A1 and A2 we obtain
\begin{eqnarray}
& & \Prob (J^{\alpha}_{N^{\alpha}(s)+1} = j, X^{\alpha}_{N^{\alpha}(s)+1} = k + v_{\alpha} \ | \ \mathbf{J}_{\mathbf{N}(s)} = \mathbf{i} \ , \ \mathbf{T}_{\mathbf{N}(s)} = s - \mathbf{v} \ , \ \mathbf{T}_{\mathbf{N}(s) + \mathbf{1}} > s )  \nonumber \\
& = & \Prob (X^{\alpha}_{N^{\alpha}(s)+1} = k + v_{\alpha} \ | \ J^{\alpha}_{N^{\alpha}(s)} = i_{\alpha} \ , \ T^{\alpha}_{N^{\alpha}(s)} = s - v_{\alpha} \ , \ T^{\alpha}_{\mathbf{N}(s) + \mathbf{1}} > s )  \nonumber \\
& \times & \Prob (J^{\alpha}_{N^{\alpha}(s)+1} = j \ | \ \mathbf{J}_{\mathbf{N}(s)} = \mathbf{i} ) ,
\end{eqnarray}
and by using the definitions and formula (\ref{Eq:HequalF}) we get
\begin{eqnarray}
& & \Prob (J^{\alpha}_{N^{\alpha}(s)+1} = j, X^{\alpha}_{N^{\alpha}(s)+1} = k + v_{\alpha} \ | \ \mathbf{J}_{\mathbf{N}(s)} = \mathbf{i} \ , \ \mathbf{T}_{\mathbf{N}(s)} = s - \mathbf{v} \ , \ \mathbf{T}_{\mathbf{N}(s) + \mathbf{1}} > s )  \nonumber \\
& = & \frac{F^{\alpha}_{i_{\alpha}}(k+v_{\alpha}) - F^{\alpha}_{i_{\alpha}}(k+v_{\alpha} - 1)}{1 - F^{\alpha}_{i_{\alpha}}(v_{\alpha})} \cdot p_{\mathbf{i},j}^{\alpha} = \frac{F^{\alpha}_{i_{\alpha}}(k+v_{\alpha}) - F^{\alpha}_{i_{\alpha}}(k+v_{\alpha} - 1)}{1 - H^{\alpha}_{i_{\alpha}}(v_{\alpha})} \cdot p_{\mathbf{i},j}^{\alpha} .
\end{eqnarray}
\end{proof}
\indent Notice that if the backward process is zero $B^{\alpha}(s) = 0$ ($v_{\alpha} = 0$) we recover
$$q^{\alpha}_{\mathbf{i},j} (0,k) = q^{\alpha}_{\mathbf{i},j} (k) .$$
Formula (\ref{Def:qBack}) reveals that assumption A1 and A2 imply that the one-step probabilities $q^{\alpha}$ are affected by the duration in the state only of the component $\alpha$. That is, the backward values of the other component does not affect the kernel $q^{\alpha}$. \\
\indent
Anyway, as we will show, the backward values of each component affect the transition probabilities of each other components. \\
\indent
Now we want to discuss the evolution equation for the MVSMC with backward recurrence times. First of all, we define the transition probability for the component $\alpha$ as
\begin{eqnarray}
\label{probab}
\Phi^{\alpha}_{\mathbf{i};j_{\alpha}} (\mathbf{v},u_{\alpha},k) := \Prob (Z^{\alpha}(k)=j_{\alpha},B^{\alpha}(k)=u_{\alpha} \ | \ \mathbf{Z}(0) = \mathbf{i} \ , \ \mathbf{B}(0)=\mathbf{v}).
\end{eqnarray}
\indent
The following result consists in a recursive system of equations for computing the transition probability functions in a bivariate system.
\begin{proposition}\label{Prop:TransEqMSMCBack}
Suppose that the system is composed of two components, i.e. $\gamma = 2$, and such that hypotheses A1, A2 hold true. Then, for all $\mathbf{i}, \mathbf{j} \in E^{2}$, $\mathbf{v}, \mathbf{u} \in \Natural^{2}$ and $k \in \Natural$, we have
\begin{equation}
\label{sistema}
\left\{ \begin{array}{l} 
\displaystyle
\Phi^{1}_{\mathbf{i};j_{1}} (\mathbf{v},u_{1},k) = \delta_{i_{1}j_{1}} \mathbb{I}_{\{u_{1}=k+v_{1}\}} \frac{1 - H^{1}_{i_{1}}(k+v_{1})}{1 - H^{1}_{i_{1}}(v_{1})} \\
\displaystyle
+ \sum_{\tau = 1}^{k} \sum_{l_{1},l_{2} \in E} \sum_{w=0}^{\tau + v_{2}} \Phi^{1}_{(l_{1},l_{2}),j_{1}} ((0,w),u_{1},k-\tau) \Phi^{2}_{(i_{1},i_{2}),l_{2}} ((v_{1},v_{2}),w,\tau) q^{1}_{(i_{1},i_{2}),l_{1}}(v_{1},\tau) \\
\displaystyle
\Phi^{2}_{\mathbf{i};j_{2}} (\mathbf{v},u_{2},k) = \delta_{i_{2}j_{2}} \mathbb{I}_{\{u_{2}=k+v_{2}\}} \frac{1 - H^{2}_{i_{2}}(k+v_{2})}{1 - H^{2}_{i_{2}}(v_{2})} \\
\displaystyle
+ \sum_{\tau = 1}^{k} \sum_{l_{1},l_{2} \in E} \sum_{w=0}^{\tau + v_{1}} \Phi^{2}_{(l_{1},l_{2}),j_{2}} ((w,0),u_{2},k-\tau) \Phi^{1}_{(i_{1},i_{2}),l_{1}} ((v_{1},v_{2}),w,\tau) q^{2}_{(i_{1},i_{2}),l_{2}}(v_{2},\tau) 
\end{array} \right.
\end{equation}
\end{proposition}
\begin{proof}
We show the result for transition function of component $1$, all what follows holds for component $2$ by symmetric arguments. \\
\indent
For all $\mathbf{i}, \mathbf{j} \in E^{2}$, $\mathbf{v}, \mathbf{u} \in \Natural^{2}$ and $k \in \Natural$, we make a partition of the sample space in two parts: no transition up to time $k$ and at least one transition; therefore we have
\begin{eqnarray}
& & \Phi^{1}_{\mathbf{i};j_{1}} (\mathbf{v},u_{1},k) = \Prob (Z^{1}(k)=j_{1},B^{1}(k)=u_{1} \ | \ \mathbf{Z}(0) = \mathbf{i} \ , \ \mathbf{B}(0)=\mathbf{v}) \nonumber \\
& = & \Prob (Z^{1}(k)=j_{1},B^{1}(k)=u_{1},T^{1}_{1}>k \ | \ \mathbf{Z}(0) = \mathbf{i} \ , \ \mathbf{B}(0)=\mathbf{v}) \nonumber \\
& + & \sum_{\tau=1}^{k} \Prob (Z^{1}(k)=j_{1},B^{1}(k)=u_{1},T^{1}_{1}=\tau \ | \ \mathbf{Z}(0) = \mathbf{i} \ , \ \mathbf{B}(0)=\mathbf{v}) .
\label{Eq:1DimPropBack}
\end{eqnarray}
\indent Let us consider the first term on the right hand side, using Bayes rules we obtain
\begin{eqnarray}
& & \Prob (Z^{1}(k)=j_{1},B^{1}(k)=u_{1},T^{1}_{1}>k \ | \ \mathbf{Z}(0) = \mathbf{i} \ , \ \mathbf{B}(0)=\mathbf{v}) \nonumber \\
& = & \Prob (Z^{1}(k)=j_{1},B^{1}(k)=u_{1} \ | \ T^{1}_{1} > k \ , \ \mathbf{Z}(0) = \mathbf{i} \ , \ \mathbf{B}(0)=\mathbf{v}) \nonumber \\
& \times & \Prob (T^{1}_{1}>k \ | \ \mathbf{Z}(0) = \mathbf{i} \ , \ \mathbf{B}(0)=\mathbf{v}) \nonumber \\
& = & \delta_{i_{1}j_{1}} \mathbb{I}_{\{u_{1}=k+v_{1}\}} \frac{1 - H^{1}_{i_{1}}(k+v_{1})}{1 - H^{1}_{i_{1}}(v_{1})} . 
\end{eqnarray}
\indent Now, let us consider the second term on the right hand side of equation (\ref{Eq:1DimPropBack}), we have
\begin{eqnarray}
& & \Prob (Z^{1}(k)=j_{1},B^{1}(k)=u_{1},T^{1}_{1}=\tau \ | \ \mathbf{Z}(0) = \mathbf{i} \ , \ \mathbf{B}(0)=\mathbf{v}) \nonumber \\
& = & \sum_{l_{1},l_{2} \in E} \sum_{w=0}^{\tau + v_{2}}  \Prob (Z^{1}(k)=j_{1},B^{1}(k)=u_{1},T^{1}_{1}=\tau,\mathbf{Z}(\tau)=\mathbf{l}, \nonumber \\
& & B^{2}(\tau)=w \ | \ \mathbf{Z}(0) = \mathbf{i} \ , \ \mathbf{B}(0)=\mathbf{v}) ,
\end{eqnarray}
where $T^{1}_{1} = \tau$ imply $B^{1}(\tau) = 0$. Then by applying Bayes rules we obtain
\begin{eqnarray}
& & \Prob (Z^{1}(k)=j_{1},B^{1}(k)=u_{1},T^{1}_{1}=\tau \ | \ \mathbf{Z}(0) = \mathbf{i} \ , \ \mathbf{B}(0)=\mathbf{v}) \nonumber \\
& = & \sum_{l_{1},l_{2} \in E} \sum_{w=0}^{\tau + v_{2}} \Prob (Z^{1}(k)=j_{1},B^{1}(k)=u_{1} \ | \ \mathbf{Z}(\tau) = \mathbf{l}, T^{1}_{1} = \tau , B^{2}(\tau) = w) \nonumber \\
& \times & \Prob( Z^{1}(\tau)=l_{1}, T^{1}_{1} = \tau \ | \ \mathbf{Z}(0) = \mathbf{i} , B^{1}(0) = v_{1}, B^{2}(\tau) = w , Z^{2}(\tau)=l_{2}) \nonumber \\
& \times & \Prob( Z^{2}(\tau)=l_{2}, B^{2}(\tau) = w \ | \ \mathbf{Z}(0) = \mathbf{i} , \mathbf{B}(0)=\mathbf{v}) .
\label{Eq:2DimBack}
\end{eqnarray}
\indent Here, by using the time homogeneity, the first term in the right hand side become
\begin{eqnarray}
\hspace{-30pt}\Prob (Z^{1}(k)=j_{1},B^{1}(k)=u_{1}|\mathbf{Z}(\tau) = \mathbf{l}, T^{1}_{1} = \tau , B^{2}(\tau) = w) = \Phi^{1}_{(l_{1},l_{2}),j_{1}} ((0,w),u_{1},k-\tau),
\label{Eq:2.1DimBack}
\end{eqnarray}
the last term, by definition, is
\begin{eqnarray}
\Prob( Z^{2}(\tau)=l_{2}, B^{2}(\tau) = w \ | \ \mathbf{Z}(0) = \mathbf{i} , \mathbf{B}(0)=\mathbf{v}) = \Phi^{2}_{(i_{1},i_{2}),l_{2}} ((v_{1},v_{2}),w,\tau) .
\label{Eq:2.2DimBack}
\end{eqnarray}
\indent Finally, let us consider the second term, we have
\begin{eqnarray}
& & \Prob( Z^{1}(\tau)=l_{1}, T^{1}_{1} = \tau \ | \ \mathbf{Z}(0) = \mathbf{i} , B^{1}(0) = v_{1}, B^{2}(\tau) = w , Z^{2}(\tau)=l_{2})  \nonumber \\
&  = & \Prob( J^{1}_{N^{1}(\tau)}=l_{1}, T^{1}_{N^{1}(\tau)} = \tau \ | \ \mathbf{Z}(0) = \mathbf{i} , B^{1}(0) = v_{1}, B^{2}(\tau) = w , Z^{2}(\tau)=l_{2})  \nonumber \\
& = & \Prob( T^{1}_{N^{1}(\tau)} = \tau \ | \ J^{1}_{N^{1}(\tau)}=l_{1}, \mathbf{Z}(0) = \mathbf{i} , B^{1}(0) = v_{1}, B^{2}(\tau) = w , Z^{2}(\tau)=l_{2}) \nonumber \\
& \times & \Prob( J^{1}_{N^{1}(\tau)}=l_{1} \ | \ \mathbf{Z}(0) = \mathbf{i} , B^{1}(0) = v_{1}, B^{2}(\tau) = w , Z^{2}(\tau)=l_{2}) .
\end{eqnarray}
\indent Here, by using assumption A2, we obtain
\begin{eqnarray}
& & \Prob( T^{1}_{N^{1}(\tau)} = \tau \ | \ J^{1}_{N^{1}(\tau)}=l_{1}, \mathbf{Z}(0) = \mathbf{i} , B^{1}(0) = v_{1}, B^{2}(\tau) = w , Z^{2}(\tau)=l_{2}) \nonumber \\
& = & \Prob( T^{1}_{N^{1}(\tau)} = \tau \ | \ J^{1}_{N^{1}(\tau)}=l_{1}, Z^{1}(0) = i_{1} , B^{1}(0) = v_{1}) \nonumber \\
& = & \Prob( T^{1}_{N^{1}(\tau)} = \tau \ | \ Z^{1}(0) = i_{1} , B^{1}(0) = v_{1}) = \frac{F^{1}_{i_{1}}(\tau + v_{\alpha}) - F^{1}_{i_{1}}(\tau + v_{1} - 1)}{1 - F^{1}_{i_{1}}(v_{1})} \nonumber \\
& = & \frac{F^{1}_{i_{1}}(\tau + v_{\alpha}) - F^{1}_{i_{1}}(\tau + v_{1} - 1)}{1 - H^{1}_{i_{1}}(v_{1})} ,
\end{eqnarray}
furthermore, by using A1, we get
\begin{eqnarray}
& & \Prob( J^{1}_{N^{1}(\tau)}=l_{1} \ | \ \mathbf{Z}(0) = \mathbf{i} , B^{1}(0) = v_{1}, B^{2}(\tau) = w , Z^{2}(\tau)=l_{2})  \nonumber \\
& = & \Prob( J^{1}_{N^{1}(\tau)}=l_{1} \ | \ \mathbf{J}_{\mathbf{N}(0)} = \mathbf{i}) = p^{1}_{\mathbf{i},l_{1}} .
\end{eqnarray}
\indent Then
\begin{eqnarray}
& & \Prob( Z^{1}(\tau)=l_{1}, T^{1}_{1} = \tau \ | \ \mathbf{Z}(0) = \mathbf{i} , B^{1}(0) = v_{1}, B^{2}(\tau) = w , Z^{2}(\tau)=l_{2}) \nonumber \\
& = & \frac{F^{1}_{i_{1}}(\tau + v_{\alpha}) - F^{1}_{i_{1}}(\tau + v_{1} - 1)}{1 - H^{1}_{i_{1}}(v_{1})} \cdot p^{1}_{\mathbf{i},l_{1}} = q^{1}_{\mathbf{i},l_{1}} (v_{1},\tau).
\label{Eq:2.3DimBack}
\end{eqnarray}
\indent By substituting the formulas (\ref{Eq:2.1DimBack},\ref{Eq:2.2DimBack},\ref{Eq:2.3DimBack}) in (\ref{Eq:2DimBack}) we get
\begin{eqnarray}
& & \Prob (Z^{1}(k)=j_{1},B^{1}(k)=u_{1},T^{1}_{1}=\tau \ | \ \mathbf{Z}(0) = \mathbf{i} \ , \ \mathbf{B}(0)=\mathbf{v}) \nonumber \\
& &\hspace{-20pt} = \sum_{l_{1},l_{2} \in E} \sum_{w=0}^{\tau + v_{2}} \Phi^{1}_{(l_{1},l_{2}),j_{1}} ((0,w),u_{1},k-\tau) \Phi^{2}_{(i_{1},i_{2}),l_{2}} ((v_{1},v_{2}),w,\tau) q^{1}_{(i_{1},i_{2}),l_{1}}(v_{1},\tau).
\end{eqnarray}
\indent Finally, by insert this last equation in (\ref{Eq:1DimPropBack}) we get the result.
\end{proof}
\noindent

\subsection{The algorithm}

To apply the model we must solve the system of equations $(\ref{sistema})$. The main steps of the algorithm are here described. As first step, by putting $k=1$ we can have only two possible values for the final backward process, that is $u_{1}=0$ or $u_{1}=v_{1}+1$.\\
\indent Let consider first the case $k=1$ and $u_{1}=0$ for the first component.\\
\indent By considering that $u_{1}=0$ implies that $\mathbb{I}_{\{u_{1}=1+v_{1}\}}=0$ we can rewrite the first equation of system $(\ref{sistema})$ as follows:
\begin{equation}
\begin{array}{l} 
\displaystyle
\Phi^{1}_{\mathbf{i};j_{1}} (\mathbf{v},0,1) = \sum_{l_{1},l_{2} \in E} \sum_{w=0}^{1 + v_{2}} q^{1}_{(i_{1},i_{2}),l_{1}}(v_{1},1) \Phi^{2}_{(i_{1},i_{2}),l_{2}} ((v_{1},v_{2}),w,1)\Phi^{1}_{(l_{1},l_{2}),j_{1}} ((0,w),0,0).
\end{array} 
\end{equation}
\indent Now, notice that $\Phi^{1}_{(l_{1},l_{2}),j_{1}} ((0,w),0,0)=\delta_{l_{1}j_{1}}$ and by substitution we get
\begin{equation}
\begin{array}{l} 
\displaystyle
\Phi^{1}_{\mathbf{i};j_{1}} (\mathbf{v},0,1) = \sum_{l_{1},l_{2} \in E}q^{1}_{(i_{1},i_{2}),l_{1}}(v_{1},1)\delta_{l_{1}j_{1}}\sum_{w=0}^{1 + v_{2}}  \Phi^{2}_{(i_{1},i_{2}),l_{2}} ((v_{1},v_{2}),w,1)\\
\displaystyle
= \sum_{l_{1} \in E}q^{1}_{(i_{1},i_{2}),l_{1}}(v_{1},1)\delta_{l_{1}j_{1}}\sum_{l_{2} \in E}\sum_{w=0}^{1 + v_{2}} \Phi^{2}_{(i_{1},i_{2}),l_{2}} ((v_{1},v_{2}),w,1)\\
\displaystyle
=\sum_{l_{1} \in E}q^{1}_{(i_{1},i_{2}),l_{1}}(v_{1},1)\delta_{l_{1}j_{1}}=q^{1}_{(i_{1},i_{2}),j_{1}}(v_{1},1).
\end{array} 
\end{equation}
\indent The second case to dealt with is $k=1$ and $u_{1}=1+v_{1}$ for the first component.\\
\indent In this evenience we have simply that
\begin{equation}
\label{Eq:k1notrans}
\begin{array}{l} 
\displaystyle
\Phi^{1}_{\mathbf{i};j_{1}} (\mathbf{v},1+v_{1},1) = \delta_{i_{1}j_{1}} \frac{1 - H^{1}_{i_{1}}(1+v_{1})}{1 - H^{1}_{i_{1}}(v_{1})}
\end{array}.
\end{equation}
\indent In fact the second addend on the right hand side of the first equation in $(\ref{sistema})$ is zero because it contains the term $\Phi^{1}_{(l_{1},l_{2});j_{1}} ((0,w),1+v_{1},0)=0$. The reason of this last equality is that, due to the fact that $k=0$, the initial and final backward values of the first component should be equal, but here we have $0\neq 1+v_{1}$ because $v_{1}$ cannot be negative.\\
\indent This first step has been executed for the first component, a symmetric argument gives similar results for the second components.\\
\indent Observe that at time $k=1$ the transition probability of one component is affected by the state of the other component but not by the duration in the state of the latter.\\
\indent As second step, by putting $k=2$ in $(\ref{sistema})$ we have the following equation:
\begin{eqnarray}
\label{k2}
& & \Phi^{1}_{\mathbf{i};j_{1}} (\mathbf{v},u_{1},2) = \delta_{i_{1}j_{1}} \mathbb{I}_{\{u_{1}=2+v_{1}\}} \frac{1 - H^{1}_{i_{1}}(2+v_{1})}{1 - H^{1}_{i_{1}}(v_{1})} \\
\displaystyle
& + & \sum_{\tau = 1}^{2} \sum_{l_{1},l_{2} \in E} \sum_{w=0}^{\tau + v_{2}} q^{1}_{(i_{1},i_{2}),l_{1}}(v_{1},\tau) \Phi^{2}_{(i_{1},i_{2}),l_{2}} ((v_{1},v_{2}),w,\tau) \Phi^{1}_{(l_{1},l_{2}),j_{1}} ((0,w),u_{1},2-\tau) . \nonumber
\end{eqnarray}
\indent When $k=2$ and the initial backward of the first component is $B^{1}(0)=v_{1}$, then $B^{1}(2)\in \{0,1,2+v_{1}\}$, consequently we have three cases to analyse: ($k=2$ and $u_{1}=2+v_{1}$, $k=2$ and $u_{1}=1$, $k=2$ and $u_{1}=0$).\\
\indent If $k=2$ and $u_{1}=2+v_{1}$ we have   
\begin{equation}
\label{k2bis}
\begin{array}{l} 
\Phi^{1}_{\mathbf{i};j_{1}} (\mathbf{v},2+v_{1},2) = \delta_{i_{1}j_{1}} \frac{1 - H^{1}_{i_{1}}(2+v_{1})}{1 - H^{1}_{i_{1}}(v_{1})} + \\
\displaystyle
+ \sum_{l_{1},l_{2} \in E} \sum_{w=0}^{1 + v_{2}} q^{1}_{(i_{1},i_{2}),l_{1}}(v_{1},1) \Phi^{2}_{(i_{1},i_{2}),l_{2}} ((v_{1},v_{2}),w,1) \Phi^{1}_{(l_{1},l_{2}),j_{1}} ((0,w),2+v_{1},1)\\
\displaystyle
+\sum_{l_{1},l_{2} \in E} \sum_{w=0}^{2 + v_{2}} q^{1}_{(i_{1},i_{2}),l_{1}}(v_{1},2) \Phi^{2}_{(i_{1},i_{2}),l_{2}} ((v_{1},v_{2}),w,2) \Phi^{1}_{(l_{1},l_{2}),j_{1}} ((0,w),2+v_{1},0) .
\end{array}
\end{equation}
\indent Note that the third added on the r.h.s. of $(\ref{k2bis})$ vanishes due to the fact that \\
\noindent
$\Phi^{1}_{(l_{1},l_{2}),j_{1}} ((0,w),2+v_{1},0)=0$ being $2+v_{1}\neq 0$. Then, It result that
\begin{equation}
\begin{array}{l} 
\Phi^{1}_{\mathbf{i};j_{1}} (\mathbf{v},2+v_{1},2) = \delta_{i_{1}j_{1}} \frac{1 - H^{1}_{i_{1}}(2+v_{1})}{1 - H^{1}_{i_{1}}(v_{1})} + \\
\displaystyle
+ \sum_{l_{1},l_{2} \in E} \sum_{w=0}^{1 + v_{2}} q^{1}_{(i_{1},i_{2}),l_{1}}(v_{1},1) \Phi^{2}_{(i_{1},i_{2}),l_{2}} ((v_{1},v_{2}),w,1) \Phi^{1}_{(l_{1},l_{2}),j_{1}} ((0,w),2+v_{1},1) ,
\end{array}
\end{equation}
moreover note that the the third factor on second addend on the right hand side vanishes, in fact by definition
\begin{eqnarray}
 \Phi^{1}_{(l_{1},l_{2}),j_{1}} ((0,w),2+v_{1},1) = \Prob (Z^{1}(1)=j_{1},B^{1}(1)=2+v^{1} \ | \ \mathbf{Z}(0) = \mathbf{l} \ , \ \mathbf{B}(0)=(0,w)) = 0
\end{eqnarray}
and finally we get
\begin{eqnarray}
\Phi^{1}_{\mathbf{i};j_{1}} (\mathbf{v},2+v_{1},2) = \delta_{i_{1}j_{1}} \frac{1 - H^{1}_{i_{1}}(2+v_{1})}{1 - H^{1}_{i_{1}}(v_{1})} ,
\end{eqnarray}
as we can expect, the probability to go in 2 unit times from an initial backward state $v_{1}$ to a final backward state $v_{1} + 2$ is the probability to make no transition in the 2 unit times. \\
\indent The second case is when $k=2$ and $u_{1}=1$ where we have 
\begin{equation}
\begin{array}{l} 
\Phi^{1}_{\mathbf{i};j_{1}} (\mathbf{v},1,2) = \delta_{i_{1}j_{1}} \mathbb{I}_{\{1=2+v_{1}\}}  \frac{1 - H^{1}_{i_{1}}(2+v_{1})}{1 - H^{1}_{i_{1}}(v_{1})} + \\
\displaystyle
+ \sum_{l_{1},l_{2} \in E} \sum_{w=0}^{1 + v_{2}} q^{1}_{(i_{1},i_{2}),l_{1}}(v_{1},1) \Phi^{2}_{(i_{1},i_{2}),l_{2}} ((v_{1},v_{2}),w,1) \Phi^{1}_{(l_{1},l_{2}),j_{1}} ((0,w),1,1)\\
\displaystyle
+\sum_{l_{1},l_{2} \in E} \sum_{w=0}^{2 + v_{2}} q^{1}_{(i_{1},i_{2}),l_{1}}(v_{1},2) \Phi^{2}_{(i_{1},i_{2}),l_{2}} ((v_{1},v_{2}),w,2) \Phi^{1}_{(l_{1},l_{2}),j_{1}} ((0,w),1,0) .
\end{array}
\end{equation}
\indent This last equation, due to the fact that $\mathbb{I}_{\{1=2+v_{1}\}}=0$ and \\
\noindent $\Phi^{1}_{(l_{1},l_{2}),j_{1}} ((0,w),1,0)=0$, simplifies in 
\begin{eqnarray}
\Phi^{1}_{\mathbf{i};j_{1}} (\mathbf{v},1,2) = \sum_{l_{1},l_{2} \in E} \sum_{w=0}^{1 + v_{2}} q^{1}_{(i_{1},i_{2}),l_{1}}(v_{1},1) \Phi^{2}_{(i_{1},i_{2}),l_{2}} ((v_{1},v_{2}),w,1) \Phi^{1}_{(l_{1},l_{2}),j_{1}} ((0,w),1,1) , \nonumber
\end{eqnarray}
now, form Step 1 (cf. Eq.(\ref{Eq:k1notrans})), the third factor on right hand side is given by
\begin{eqnarray}
\Phi^{1}_{(l_{1},l_{2}),j_{1}}((0,w),1,1) = \delta_{l_{1},j_{1}} (1 - H^{1}_{l_{1}}(1)) = \delta_{l_{1},j_{1}} (1 - H^{1}_{j_{1}}(1)) ,
\end{eqnarray}
then, by substitution, we finally get
\begin{eqnarray}
\Phi^{1}_{\mathbf{i};j_{1}} (\mathbf{v},1,2) =  q^{1}_{(i_{1},i_{2}),j_{1}}(v_{1},1) (1 - H^{1}_{j_{1}}(1)) ,
\end{eqnarray}
as we can expect, the probability to go in 2 unit times from an initial backward state $v_{1}$ to a final backward state $1$ is just the probability to make one transition on the first unit time and then wait one  uint time. \\
\indent The last case is $k=2$ and $u_{1}=0$ where we have 
\begin{equation}
\label{k2tris}
\begin{array}{l} 
\Phi^{1}_{\mathbf{i};j_{1}} (\mathbf{v},0,2) = \delta_{i_{1}j_{1}} \mathbb{I}_{\{0=2+v_{1}\}}  \frac{1 - H^{1}_{i_{1}}(2+v_{1})}{1 - H^{1}_{i_{1}}(v_{1})} + \\
\displaystyle
+ \sum_{l_{1},l_{2} \in E} \sum_{w=0}^{1 + v_{2}} q^{1}_{(i_{1},i_{2}),l_{1}}(v_{1},1) \Phi^{2}_{(i_{1},i_{2}),l_{2}} ((v_{1},v_{2}),w,1) \Phi^{1}_{(l_{1},l_{2}),j_{1}} ((0,w),0,1)\\
\displaystyle
+\sum_{l_{1},l_{2} \in E} \sum_{w=0}^{2 + v_{2}} q^{1}_{(i_{1},i_{2}),l_{1}}(v_{1},2) \Phi^{2}_{(i_{1},i_{2}),l_{2}} ((v_{1},v_{2}),w,2) \Phi^{1}_{(l_{1},l_{2}),j_{1}} ((0,w),0,0) .
\end{array}
\end{equation}
\indent Remark that $\mathbb{I}_{\{0=2+v_{1}\}}=0$ and $\Phi^{1}_{(l_{1},l_{2}),j_{1}} ((0,w),0,0)=\delta_{l_{1}j_{1}}$ make $(\ref{k2tris})$ equal to
\begin{eqnarray}
& & \Phi^{1}_{\mathbf{i};j_{1}} (\mathbf{v},0,2) = 0+ \sum_{l_{1},l_{2} \in E} \sum_{w=0}^{1 + v_{2}} q^{1}_{(i_{1},i_{2}),l_{1}}(v_{1},1) \Phi^{2}_{(i_{1},i_{2}),l_{2}} ((v_{1},v_{2}),w,1) \Phi^{1}_{(l_{1},l_{2}),j_{1}} ((0,w),0,1) \nonumber \\
\displaystyle
& & +\sum_{l_{1},l_{2} \in E} \sum_{w=0}^{2 + v_{2}} q^{1}_{(i_{1},i_{2}),l_{1}}(v_{1},2) \Phi^{2}_{(i_{1},i_{2}),l_{2}} ((v_{1},v_{2}),w,2) \delta_{(l_{1},j_{1})} \nonumber \\
\displaystyle
& & =\sum_{l_{1},l_{2} \in E} \sum_{w=0}^{1 + v_{2}} q^{1}_{(i_{1},i_{2}),l_{1}}(v_{1},1) \Phi^{2}_{(i_{1},i_{2}),l_{2}} ((v_{1},v_{2}),w,1) \Phi^{1}_{(l_{1},l_{2}),j_{1}} ((0,w),0,1) \\
\displaystyle
& & +q^{1}_{(i_{1},i_{2}),j_{1}}(v_{1},2) \sum_{l_{2} \in E} \sum_{w=0}^{2 + v_{2}} \Phi^{2}_{(i_{1},i_{2}),l_{2}} ((v_{1},v_{2}),w,2)\nonumber \\
\displaystyle
& & \hspace{-10pt} =\sum_{l_{1},l_{2} \in E} \sum_{w=0}^{1 + v_{2}} q^{1}_{(i_{1},i_{2}),l_{1}}(v_{1},1) \Phi^{2}_{(i_{1},i_{2}),l_{2}} ((v_{1},v_{2}),w,1) \Phi^{1}_{(l_{1},l_{2}),j_{1}} ((0,w),0,1) + q^{1}_{(i_{1},i_{2}),j_{1}}(v_{1},2) , \nonumber
\end{eqnarray}
which are all known terms from Step 1. A symmetric argument gives similar results for the second components. At this point, proceeding this forward algorithm we can get $\Phi^{\alpha}_{\mathbf{i};j_{1}} (\mathbf{v},u_{\alpha},k)$ by knowing $\Phi^{\beta}_{\mathbf{i};j_{1}} (\mathbf{v},u_{\beta},s)$ for all $s<k$ and $\beta \in \{1,2\}$.

\section{The 2-component reliability model}
 
Let us consider two partitions of the state space $E$ of the two components as follows:
\begin{equation}
\begin{array}{l} 
E=U^{1}\bigcup D^{1},\,\,U^{1}\bigcap D^{1}=\emptyset;\\
\displaystyle
E=U^{2}\bigcup D^{2},\,\,U^{2}\bigcap D^{2}=\emptyset.\\
\end{array},
\end{equation}
The subset $U^{1}$ ($U^{2}$) contains all good states in which the component $1$ ($2$) is regarded as well working and the subset $D^{1}$ ($D^{2}$) contains all the bad states in which the first (second) component is not well performing.\\
\indent Let us formulate an additional assumption named A3.
\begin{itemize}
\item[A3] All states in $D^{1}$ ($D^{2}$) are absorbing for the component $1$ ($2$).
\end{itemize}
\indent Notice that, the assumption A3 can be relaxed with easiness and it is possible to execute the following computation in the general case in which the first (second) component alternates between the set $U^{1}$ and $D^{1}$ ($U^{2}$ and $D^{2}$). Anyway we adopt A3 because the presentation is easier and furthermore the application we discuss in next section falls well within this case.\\ 
\indent One of the most useful indicators is the reliability function. Here we define the reliability of the system as follows:
\begin{equation}
\label{relia2}
R_{i_{1}i_{2}}((v_{1},v_{1}),k):= \Prob (Z^{1}(k)\in U^{1},Z^{2}(k)\in U^{2} \ | \ \mathbf{Z}(0)=\mathbf{i}, \mathbf{B}(0)=\mathbf{v}). 
\end{equation}
\indent In what follows we need also of the marginal reliability functions for a single component defined as:
\begin{equation}
\label{relia}
\begin{array}{l}
R_{i_{1},i_{2}}^{1}((v_{1},v_{2}),k):= \Prob (Z^{1}(k)\in U^{1} \ | \ \mathbf{Z}(0)=\mathbf{i}, \mathbf{B}(0)=\mathbf{v})\\
\displaystyle
R_{i_{1},i_{2}}^{2}((v_{1},v_{2}),k):= \Prob (Z^{2}(k)\in U^{2} \ | \ \mathbf{Z}(0)=\mathbf{i}, \mathbf{B}(0)=\mathbf{v}) .
\end{array}
\end{equation}
\indent The following results expresses the link between the marginal reliabilities and the system reliability in our model:
\begin{proposition}
\label{Prop:Reliability}
Suppose that the system is composed of two components, i.e. $\gamma = 2$, and such that hypotheses A1, A2 and A3 hold true. Then, for all $\mathbf{i}, \mathbf{j} \in E^{2}$, $\mathbf{v}, \in \Natural^{2}$ and $k \in \Natural$, we have
\begin{equation}
\label{relibilitycalcolo}
R_{i_{1}i_{2}}((v_{1},v_{2}),k)=R_{i_{1},i_{2}}^{1}((v_{1},v_{2}),k)R_{i_{1},i_{2}}^{2}((v_{1},v_{2}),k)
\end{equation}
\end{proposition}
\noindent where for $\alpha \in \{1,2\}$ we have
\begin{eqnarray}
\displaystyle
 R_{i_{1},i_{2}}^{\alpha}((v_{1},v_{2}),k)= \sum_{u_{\alpha} \geq 0}\sum_{j_{\alpha}\in U^{\alpha}}\Phi^{\alpha}_{\mathbf{i};j_{\alpha}} (\mathbf{v},u_{\alpha},k) = \sum_{j_{\alpha}\in U^{\alpha}}\Phi^{\alpha}_{\mathbf{i};j_{\alpha}} (\mathbf{v},\cdot,k) .
\end{eqnarray}
\begin{proof}
By using Bayes rule we get
\begin{equation}
\begin{array}{l}
R_{i_{1}i_{2}}((v_{1},v_{1}),k)= \Prob (Z^{1}(k)\in U^{1},Z^{2}(k)\in U^{2} \ | \ \mathbf{Z}(0)=\mathbf{i}, \mathbf{B}(0)=\mathbf{v})\\
\displaystyle
= \Prob (Z^{1}(k)\in U^{1} \ | \ Z^{2}(k)\in U^{2}, \mathbf{Z}(0)=\mathbf{i}, \mathbf{B}(0)=\mathbf{v}) \Prob (Z^{2}(k)\in U^{2} \ | \ \mathbf{Z}(0)=\mathbf{i}, \mathbf{B}(0)=\mathbf{v}) ,
\end{array}
\end{equation}
and by using formula (\ref{probab}), the definitions of the semi-Markov and backward processes and assumption A3, we get in 
\begin{equation}
\begin{array}{l}
R_{i_{1}i_{2}}((v_{1},v_{1}),k)=\sum_{j_{2}\in U^{2}}\Phi^{2}_{\mathbf{i};j_{2}} (\mathbf{v},\cdot,k)\\
\displaystyle
\times \Prob (J_{N^{1}(k)}^{1}\in U^{1} \ | \ J_{N^{2}(k)}^{2}\in U^{2}, J_{N^{1}(0)}^{1}=i_{1}, J_{N^{2}(0)}^{2}=i_{2}, T_{N^{1}(0)}^{1}=-v_{1},T_{N^{2}(0)}^{2}=-v_{2}) .
\end{array}
\end{equation}
\indent By using assumption A1 and A2 we have
\begin{equation}
\begin{array}{l}
R_{i_{1}i_{2}}((v_{1},v_{1}),k)\\
\displaystyle
= \Prob (J_{N^{1}(k)}^{1}\in U^{1} \ | \  J_{N^{1}(0)}^{1}=i_{1}, J_{N^{2}(0)}^{2}=i_{2}, T_{N^{1}(0)}^{1}=-v_{1},T_{N^{2}(0)}^{2}=-v_{2})\sum_{j_{2}\in U^{2}}\Phi^{2}_{\mathbf{i};j_{2}} (\mathbf{v},\cdot,k)\\
\displaystyle
=\sum_{j_{1}\in U^{1}}\Phi^{1}_{\mathbf{i};j_{1}} (\mathbf{v},\cdot,k)\sum_{j_{2}\in U^{2}}\Phi^{2}_{\mathbf{i};j_{2}} (\mathbf{v},\cdot,k)=R_{i_{1}}^{1}((v_{1},v_{2}),k)R_{i_{1}}^{2}((v_{1},v_{2}),k) .
\end{array}
\end{equation}
\end{proof}
\begin{rem}
The precedent result allow us to express all the reliabilities in terms of the transition probabilities.
\label{Rem:Rel1}
\end{rem}
\begin{rem}
If we consider a model where the two components are supposed to be independent, then, the reliability of the single component would be modelled with a standard univariate semi-Markov chain and we would get
\begin{equation}
\label{relia}
\begin{array}{l}
R_{i_{1}}^{1}(v_{1},k):= \Prob (Z^{1}(k)\in U^{1} \ | \ Z^{1}(0)=i_{1}, B^{1}(0)=v_{1})\\
\displaystyle
R_{i_{2}}^{1}(v_{2},k):= \Prob (Z^{2}(k)\in U^{2} \ | \ Z^{2}(0)=i_{2}, B^{2}(0)=v_{2}) .
\end{array}
\end{equation}
If the two components are supposed to be independent, the product of these two reliabilities should be equal to the product of the bivariate reliabilities evaluated in Proposition \ref{Prop:Reliability}.
Then, any deviation of the ratio 
$$\frac{R_{i_{1}i_{2}}^{1}((v_{1},v_{2}),k)R_{i_{1}i_{2}}^{2}((v_{1},v_{2}),k)}{R_{i_{1}}^{1}(v_{1},k)R_{i_{2}}^{2}(v_{2},k)}$$
 by one is an indication of the correlation between the two components. 
 \label{Rem:Rel2}
 \end{rem}

\section{Counterparty Credit Risk in a CDS contract} 
\label{Sec:CCRforCDS}

In the financial market all subjects are exposed to the default risk. Then, in any financial contract we have to take into account for the risk of default of the our counterpart. Counterpary credit risk is in general `the risk that a counterpart of a financial contract will default prior to the expiration of the contract and will not make all the payments required by the contract' (cf. \textsl{Pykhtin and Zhu} \cite{Zhu}). \\
\noindent
We would like to study the counterparty credit risk in a Credit Default Swap (CDS) contract. In this work we would like to emphasize the difference between the CDS contract with and without consider counterparty risk, we will call these two cases \textsl{risky CDS} and \textsl{risk free CDS} respectively (\textsl{Cr\'{e}pey et al.} \cite{Crepey}). \\
\noindent
Let us consider a firm C, supposed to be defaultable, emitting an obligation (or bond) on one money unit at the time $0$ with maturity time $T$. Let us also consider a bondholder A (or protection buyer) supposed to be risk free in all what follows. The possible financial scenarios are
\begin{itemize}
\item If C has not been default until T, it is able to pay the money due to the bondholder A.
\item In case of C's default before or at the maturity date $T$, it will be able only to pay a fraction (recovery rate $\rho_{C}$) of the face value of the obligation to A.
\end{itemize}
For these reasons the bondholder A is looking for protection against the loss in case of C's default. Let us consider a third financial subject that we will call generically as protection seller B. A risk free CDS is a contract which obligates A (protection buyer) to pay a fee to B (protection seller supposed to be risk free) in change of protection against the default of the reference credit firm C. The cash flows of a risk free CDS are
\begin{itemize}
\item A pays to B a stream of premia with spread k, from the initial date until the occurence of default event or the maturity date $T$.
\item In case of default of C, B has to cover the loss of A. Then B has to pay $1 - \rho_{C}$ unit of money to A. 
\end{itemize}
The value of the spread is evaluated in order to guarantee, that the contract has value zero at the inception time. We assume that the payment of B to A is made at the same time of the default event. \\
\noindent
Let $\tau_{C}$ be the time of default for the credit reference firm C. From the above discussion, we can directly write an expression for the cash flows and price process of the risk free CDS contract. The In Cash Flows process from the perspective of the bondholder A in the risk free CDS is given by
\begin{eqnarray}
\beta_{\tau_{C}} (1 - \rho_{C}) \mathbb{I}_{\{t<\tau_{C}\leq T\}} ,
\end{eqnarray}
where $\beta$ is a discount factor. The Out Cash Flows process is given by
\begin{eqnarray}
- K \sum_{s = t} ^{T} \beta_{s} \mathbb{I}_{\{s<\tau_{C}\}} .
\end{eqnarray}
Then, the discounted value of the risk free CDS with maturity $T$ is
\begin{eqnarray}
\beta_{t}p_{T}(t) = - K \sum_{s = t} ^{T} \beta_{s} \mathbb{I}_{\{s<\tau_{C}\}} + \beta_{\tau_{C}} (1 - \rho_{C}) \mathbb{I}_{\{t<\tau_{C}\leq T\}} \nonumber
\end{eqnarray}
and its price process is given by $P_{t} = \mathbb{E}_{t} [p_{T}(t)]$. The subscript $t$, here and after, indicates the information that at time $t$ the process is still in one of the Up states.
\begin{rem}
The price of a risk free CDS can be evaluated with a single component reliability model where the rating of the only defaultable subject is modeled via a standard univariate semi-Markov process (see for example \textsl{D'Amico et al.} \cite{DamicoManca3}).
\label{Rem:CCR1}
\end{rem}
A risky CDS is a contract which obligates A (protection buyer) to pay a fee to B (defaultable protection seller) in change of protection against the default of the reference credit firm C. The cash flows of a risky CDS are
\begin{itemize}
\item A pays to B a stream of premia with spread k, from the initial date until the occurence of default event or the maturity date $T$.
\item In case of default of C, if B has not defaulted, B has to cover the loss of A. Then B has to pay $1 - \rho_{C}$ unit of money to A.
\item In case of default of B, if C has not defaulted, the contract is stopped with a \emph{Close-Out Cash Flow} (cf. \textsl{Cr\'{e}pey et al.} \cite{Crepey}). In this work we assume that the two parties according on a termination of the contract with a terminal cash flow paid to A, positive or negative, depending on the value of the risk free CDS computed at the time of default (cf. \textsl{Brigo et al.} \cite{Brigo}). 
\item If B defaults at the same time as the firm C, B will be only able to pay to A a fraction (recovery rate $\rho_{B}$)of the loss of A, namely $\rho_{B}(1 - \rho_{C})$ unit of money.
\end{itemize}
The value of the spread is evaluated in order to guarantee, that the contract has value zero at the inception time. We assume that the payment of B to A is made at the same time of the default event(s). 
The possible loss of A for the joint default event is an effect due to the counterparty risk. \\
\indent
Let us introduce $\tau_{B}$, the time of default for the protection seller B. The In Cash Flows process for the risky CDS is given by
\begin{eqnarray}
\beta_{\tau_{C}} (1 - \rho_{C}) \mathbb{I}_{\{t<\tau_{C}\leq T\}} [ \mathbb{I}_{\{\tau_{C} < \tau_{B} \}} + \rho_{B} \mathbb{I}_{\{\tau_{C} = \tau_{B} \}}] + \beta_{\tau_{B}} \mathbb{I}_{\{t<\tau_{B}\leq (T \wedge \tau_{C})\}} \rho_{B} P_{\tau_{B}}^{+}
\end{eqnarray}
here $\beta$ is a discount factor and with $P_{\cdot}^{+}$ we denote the positive part of the price process for the risk free CDS. The Out Cash Flows process is given by
\begin{eqnarray}
- K \sum_{s = t} ^{T} \beta_{s} \mathbb{I}_{\{s<(\tau_{C} \wedge \tau_{B})\}} - \beta_{\tau_{B}} \mathbb{I}_{\{t<\tau_{B}\leq (T \wedge \tau_{C})\}} P_{\tau_{B}}^{-}
\end{eqnarray} 
where $P_{\cdot}^{-}$ stands for the negative part of the price process for the risk free CDS. Then, the discounted value of the risky CDS with maturity $T$ is
\begin{eqnarray}
\beta_{t}\pi_{T}(t) & = & - K \sum_{s = t} ^{T} \beta_{s} \mathbb{I}_{\{s<(\tau_{C} \wedge \tau_{B})\}} + \beta_{\tau_{C}} (1 - \rho_{C}) \mathbb{I}_{\{t<\tau_{C}\leq T\}} [ \mathbb{I}_{\{\tau_{C} < \tau_{B} \}} + \nonumber \\
& + & \rho_{B} \mathbb{I}_{\{\tau_{C} = \tau_{B} \}}] + \beta_{\tau_{B}} \mathbb{I}_{\{t<\tau_{B}\leq (T \wedge \tau_{C})\}} (\rho_{B} P_{\tau_{B}}^{+} - P_{\tau_{B}}^{-}) ,
\label{Eq:RiskyPriceProc}
\end{eqnarray}
and the price process for the risky CDS is $\Pi_{t} = \mathbb{E}_{t} [\pi_{T}(t)]$.

\subsection{Pricing a risky CDS and CVA evaluation}
In this section we are going to apply the 2-component reliability model to the study of the counterpart risk in a CDS contract. In particular our goal is to price a risky CDS contract and to derive an expression for the credit value adjustment (CVA) which can be seen as a measure of the counterparty credit risk. \\
\noindent
In order to price a risky CDS we should be able to evaluate the expected value of the indicator functions in (\ref{Eq:RiskyPriceProc}). The following result concern the evaluation of the expectation.
\begin{proposition}
The price of a risky CDS under the natural probability measure is
\begin{eqnarray}
\begin{aligned}
\beta_{t}\Pi_{T}(t) & =  - K \sum_{s = t}^{T} \beta_{s} \sum_{h=s+1}^{T} \Prob_{t}(\tau = h) + \sum_{h_{B}=t+1}^{\infty} \sum_{h_{C}=t+1}^{x_{B}} \beta_{h_{C}} (1 - \rho_{C}) \Prob_{t} ( \tau_{C} = h_{C} , \tau_{B} = h_{B} ) \\
& + \sum_{h_{C}=t+1}^{T} \beta_{h_{C}} (1 - \rho_{C}) \rho_{B} \Prob_{t} ( \tau_{C} = h_{C} , \tau_{B} = h_{C} ) ,
\end{aligned}
\end{eqnarray}
where $\tau = \tau_{C} \wedge \tau_{B}$,
\begin{eqnarray}
\Prob_{t} (\tau = h) = \sum_{h_{B}=h}^{\infty} \Prob_{t} (\tau_{C} = h, \tau_{B} = h_{B}) + \sum_{h_{C}=h}^{\infty} \Prob_{t} (\tau_{C} = h_{C}, \tau_{B} = h) ,
\end{eqnarray}
and
\begin{eqnarray}
\Prob_{t} ( \tau_{C} = h_{C} , \tau_{B} = h_{B} ) = \sum_{\mathbf{i} \in U^{2}} \sum_{\mathbf{v} \in \Natural^{2}} (R^{C}_{\mathbf{i},\mathbf{v}}(h_{C} - 1) - R^{C}_{\mathbf{i},\mathbf{v}}(h_{C})) (R^{B}_{\mathbf{i},\mathbf{v}}(h_{B} - 1) - R^{B}_{\mathbf{i},\mathbf{v}}(h_{B})) .
\end{eqnarray}
\end{proposition}
\begin{proof}
The result is a direct consequence of formula (\ref{Eq:RiskyPriceProc}) but it remain to prove that
\begin{eqnarray}
\Prob_{t} ( \tau_{C} = h_{C} , \tau_{B} = h_{B} ) = \sum_{\mathbf{i} \in U^{2}} \sum_{\mathbf{v} \in \Natural^{2}} (R^{C}_{\mathbf{i},\mathbf{v}}(h_{C} - 1) - R^{C}_{\mathbf{i},\mathbf{v}}(h_{C})) (R^{B}_{\mathbf{i},\mathbf{v}}(h_{B} - 1) - R^{B}_{\mathbf{i},\mathbf{v}}(h_{B})) .
\end{eqnarray}
We will show the result for $t = 0$ for the sake of simplicity, the general case is a direct consequence. First of all we notice that
\begin{eqnarray}
\Prob_{0} ( \tau_{C} = h_{C} , \tau_{B} = h_{B} ) = \sum_{\mathbf{i} \in U^{2}} \sum_{\mathbf{v} \in \Natural^{2}} \Prob (Z^{C}(h_{C}) = D,Z^{B}(h_{B}) = D \ | \ \mathbf{Z}(0)=\mathbf{i}, \mathbf{B}(0)=\mathbf{v}).
\end{eqnarray}
Using the properties of the probability measures, we get
\begin{eqnarray}
\begin{aligned}
\Prob &(Z^{C}(h_{C}) = D,Z^{B}(h_{B}) = D \ | \ \mathbf{Z}(0)=\mathbf{i}, \mathbf{B}(0)=\mathbf{v}) = 1 - \left[ \Prob (Z^{C}(h_{C}) \in U \ | \ \mathbf{Z}(0)=\mathbf{i}, \mathbf{B}(0)=\mathbf{v}) \right. \\
&  \left. + \Prob (Z^{B}(h_{B}) \in U \ | \ \mathbf{Z}(0)=\mathbf{i}, \mathbf{B}(0)=\mathbf{v}) - \Prob (Z^{C}(h_{C}) \in U,Z^{B}(h_{B}) \in U \ | \ \mathbf{Z}(0)=\mathbf{i}, \mathbf{B}(0)=\mathbf{v}) \right] .
\end{aligned}
\end{eqnarray}
The last expression can be rewritten in terms of reliabilities as
\begin{eqnarray}
\begin{aligned}
\Prob &(Z^{C}(h_{C}) = D,Z^{B}(h_{B}) = D \ | \ \mathbf{Z}(0)=\mathbf{i}, \mathbf{B}(0)=\mathbf{v}) = (R^{C}_{\mathbf{i},\mathbf{v}}(h_{C} - 1) - R^{C}_{\mathbf{i},\mathbf{v}}(h_{C})) (R^{B}_{\mathbf{i},\mathbf{v}}(h_{B} - 1) - R^{B}_{\mathbf{i},\mathbf{v}}(h_{B})) ,
\end{aligned}
\end{eqnarray}
which concludes the proof.
\end{proof}
\begin{rem}
The precedent result gives an expression of the price of a risky CDS as a function of the reliabilities, that is the transition probabilities for the bivariate semi-Markov chain (cf. remark \ref{Rem:Rel1}).
\end{rem}
\begin{rem}
The difference between the price of a risky CDS and the price of a risk free CDS is has a particular importance, indeed it is a measure of the loss of value a CDS contract undergoes due to the counterpart credit risk. This difference is called credit value adjustment (CVA). The credit value adjustment process $(CVA_{t})$ is defined by
\begin{eqnarray}
CVA_{t} = P_{t} - \Pi_{t} \qquad \textrm{for} \  t < \tau_{B} ,
\end{eqnarray}
it measures the loss of value of the CDS contract. We notice that, in our model, the CVA process can be totally expressed in term of the reliabilities.
\end{rem}

\section{Conclusions} \label{conc}

This paper proposes a multivariate semi-Markov chain model in discrete time. The multicomponent system is analysed in the transient case by 
giving methods for computing the transition probabilities and reliability functions. The numerical solution of the equations is made available by
means of a recursive algorithm.\\ 
\indent The results are applied to the evaluation of the risky credit default swap contracts and they allow the attainment of an explicit formula for the price of a risky CDS and for the credit value adjustment process.\\
\indent Possible avenues for future developments of our model include:\\
a) application to real data on credit rating dynamics;\\
b) asymptotic properties of the multivariate semi-Markov model;\\
c) construction of a multivariate reward model for the credit spread computation.

\end{document}